\documentstyle[amscd,amssymb,12pt]{amsart}
\input xypic \xyoption{curve}

\def\hcorrection#1{\advance\hoffset by #1 }
\def\vcorrection#1{\advance\voffset by #1 }

\vcorrection{-.7in}
\hcorrection{-1in}

\newcommand{\C}[1]{{\cal#1}} 
\newcommand{\D}[1]{{\Bbb#1}} 

\theoremstyle{plain}
\newtheorem{th}{Theorem}[section]
\newtheorem{cor}{Corollary}[section]
\newtheorem{lem}{Lemma}[section]
\newtheorem{prop}{Proposition}[section]

\theoremstyle{definition}
\newtheorem{defin}{Definition}[section]

\theoremstyle{definition}

\theoremstyle{remark}
\newtheorem{rem}{Remark}[section]
\newtheorem{notation}{Notation}[section]

\numberwithin{equation}{section}

\begin{document}

\pagestyle{plain}
\addtolength{\footskip}{.3in}

\title[Graph cohomology and perturbative QFT]
{Cohomology of Feynman Graphs\\
and\\ Perturbative Quantum Field Theory}

\author{Lucian M. Ionescu}
\address{Department of Mathematics, Illinois State University, IL 61790-4520}
\email{lmiones@@ilstu.edu}
\date{12/23/2003}

\begin{abstract}
An analog of Kreimer's coproduct from renormalization of Feynman integrals
in quantum field theory,
endows an analog of Kontsevich's graph complex with a dg-coalgebra structure.
The graph complex is generated by orientation classes of labeled directed graphs.
A graded commutative product is also defined,
compatible with the coproduct.
Moreover, a dg-Hopf algebra is identified.

Graph cohomology is defined
applying the cobar construction to the dg-coalgebra structure.

As an application, 
L-infinity morphisms represented as series over Feynman graphs
correspond to graph cocycles. 
Notably the total differential of the cobar construction
corresponds to the L-infinity morphism condition.
The main example considered is Kontsevich's formality morphism.

The relation with perturbative quantum field theory is considered by
interpreting L-infinity morphisms as partition functions,
and the coefficients of the graph expansions as Feynman integrals.
\end{abstract}

\maketitle
\tableofcontents


\section{Introduction}\label{S:intro}
In this article we put together 
the Kreimer's coproduct (main organizational principle in renormalization \cite{CKren1,K2}),
with Kontsevich graph homology \cite{Kon2}, 
obtaining a dg-Hopf algebra structure on the k-module of 
orientation classes of directed labeled graphs.

The cobar construction applied to the dg-coalgebra structure yields graph cohomology,
and notably the cobar total differential corresponds to the L-infinity morphisms
expanded over the corresponding class of graphs.
As a consequence,
L-infinity morphisms correspond to graph cocycles,
Kontsevich's formality morphism being a prototypical example.

This mathematical framework is adequate for perturbative quantum field theory,
with partition functions corresponding to L-infinity morphisms and
with cocycles playing the role of Feynman integrals.

Section \ref{S:hag} introduces the algebra of graphs and defines the dg-coalgebra structure.
The additional structure yielding the dg-Hopf algebra structure is mentioned.

In section \ref{S:gc} graph cohomology is defined.

Section \ref{S:Kfm} is concerned with L-infinity morphisms as perturbation series
over a given class of graphs.
The main example of Kontsevich's formality morphism is revisited.
Prompted by the interpretation of an L-infinity morphism as a partition function,
``Feynman rules/integrals'' are defined via 
a pairing between graphs and states on graphs.

The article concludes in section \ref{S:relpQFT} 
with some general comments regarding the relation between
L-infinity algebras and perturbative QFT.

\section{The dg-Hopf algebra of oriented graphs}\label{S:hag}
Let $\C{G}$ be a class of directed graphs with orientation modulo equivalence (defined below).
Typical examples are $\{G_{n,m}\}_{n,m\ge 0}$, 
the class of admissible graphs of Kontsevich (\cite{Kon1}, p.22),
and the Feynman graphs corresponding to a given (perturbative ) QFT 
(e.g. $\phi^3_D$ \cite{CKren1}, p.9).

Denote by $\Gamma^{(1)}$ the set of (directed) {\em edges} of $\Gamma$, 
and by $\Gamma^{(0)}$ the set of {\em vertices}.
The number of {\em internal vertices}  will be denoted by $n=|\Gamma^{(0)}_{int}|$ 
and $m=|\Gamma^{(0)}_{bd}|$ will denote the number of {\em boundary vertices} 
($\partial \Gamma=\Gamma^{(0)}_{bd}$ for short).

The {\em orientation class} $[\Gamma,l]$ of a graph $\Gamma$,
is defined with respect to a {\em labeling} $l$ of its internal vertices
$l^{(0)}:\Gamma^{0}_{int}\to \{1,...,n\}$, 
and an enumeration $l^{(1)}$ of outgoing edges \cite{Kon1}, p.22,
while the labeling of boundary vertices $l^{(0)}_{bd}:\partial\Gamma\to \{1,...,m\}$
will play the role of an order.
The labeling of the boundary points will not be reflected in the notation used,
yet assumed throughout.

Two {\em labeled graphs} (as above),
obtained by a transposition of adjacent edges with the same source
represent opposite orientation classes.
When two adjacent vertices are transposed, 
a corresponding sign $(-1)^{a b}$ must be included as a coefficient of the orientation class,
where $a$ and $b$ denote the number of outgoing edges 
from the two vertices ($\#star$ in \cite{Kon1}).
\begin{rem}\label{R:orient}
The present orientation is needed to pair internal vertices with 
elements of the exterior algebra of polyvector fields (see Section \ref{S:gexp}),
while the boundary vertices will be paired with elements of a tensor algebra.

For additional details regarding orientation,
compare with the orientation on graphs which are not directed,
defined in \cite{Kon2} p.108, \cite{Kon3} p.175, \cite{CV}, p.2,
and suited for the correspondence with the punctured Riemann surfaces 
(see also \cite{Vor1} Lecture 4, p.1.)
\end{rem}
%
Let $H=k<\C{G}>$ be the graded k-module of oriented graphs,
graded by the number of edges $deg(\Gamma)=|\Gamma^{(1)}|$.

Introduce a {\em product} determined by the disjoint union of the underlying graphs and
concatenation of orders (labelings) (adapting \cite{CV}, p.9):
$$[\Gamma_1,l_1][\Gamma_2,l_2]=[\Gamma_1\cup\Gamma_2,l_1<l_2].$$
\begin{lem}
$(H, \cdot)$ is a graded commutative algebra.
\end{lem}
\begin{pf}
We will only note that (\cite{CV}, p.9):
$$[\Gamma_1,l_1] [\Gamma_2,l_2]=(-1)^{deg(\Gamma_1) deg(\Gamma_2)}
[\Gamma_2,l_2] [\Gamma_1,l_1]$$.
\end{pf}
The {\em unit} is the k-linear map $\epsilon:k\to H$ assigning to the unit of the
ground field the empty graph: $\epsilon(1)=\emptyset\in G_{0,0}$.
The empty graph will be identified as the unit of $H$, and denoted as $1$.
It is immediate to verify that it is a unit for the above defined product.

The {\em coproduct} is defined by:
\begin{equation}\label{E:coprod}
\Delta[\Gamma,l]=[\Gamma,l]\otimes 1+1\otimes[\Gamma,l]+\Delta_b[\Gamma,l],
\end{equation}
where the reduced coproduct $\Delta_b$ is defined by:
\begin{equation}\label{E:fgcoprod}
\Delta_b[\Gamma,l]=\sum_{\gamma\to \Gamma\to \gamma', \ \gamma\cap\partial\Gamma\ne \emptyset}
[\gamma,l|\gamma]\otimes[\gamma',l/\gamma].
\end{equation}
The above sum is over all {\em normal subgraphs} of $\Gamma$ meeting the boundary $\partial\Gamma$,
i.e. such that collapsing $\gamma$ to a (boundary) vertex yields a graph from the given class $\C{G}$ 
(compare being ``normal'' with condition (7) \cite{CKren1}, p.11).
Note that the additional two terms in Equation \ref{E:coprod} may be obtained 
as terms in Equation \ref{E:fgcoprod},
if allowing the subgraph $\gamma$ to be the graph itself and the empty subgraph.

%
The {\em counit} is the k-linear map $\eta:H\to k$ determined by $\eta(\Gamma)=0$, $\Gamma\ne \emptyset$
(projection onto the ``-1'' degree component of $H$, generated by the unit).
It is immediate to verify that it is a counit for the above defined coproduct.
%
\begin{rem}\label{R:whyint}
As defined above, 
the terms of the graph homology differential correspond to internal edges only,
while those of the reduced coproduct correspond only to subgraphs collapsing to a boundary vertex.
These requirements are specific of the application aiming to construct L-infinity morphisms.

It should be noted that in \cite{Kon1} the integral over the codimension one boundary of the
compactification of the configuration space of a given graph is a sum over all
its proper subgraphs (\cite{Kon1}, 5.2.1., p.22; see \cite{pqft} for additional details).
The other terms vanish for various reasons,
and the sum reduces to the above terms corresponding to $d$ and $\Delta_b$.
These in turn correspond to the L-infinity condition,
as proved in section \ref{S:Kfm}.
\end{rem}
The orientation class of the subgraph $l|\gamma$,
is defined by the restriction of a total order $l$ on $\Gamma$
corresponding to the given orientation class
(compare with \cite{CV}, p.3, \cite{Vor1}, L4, p.1, L9, p.3).

The orientation class of the collapsed graph is defined 
by the labeling induced by restriction,
except for the label of the vertex which is obtained from collapsing $\gamma$.
By changing the representative if necessary,
assume that the vertices of $\gamma$ precede the other vertices in $\Gamma$.
Then the collapsed subgraph will be the first vertex of the quotient graph $\Gamma/\gamma$.
This {\em quotient labeling} will be denoted by $l/\gamma$.
\begin{lem}
$\Delta:H\to H^{\otimes^2}$ is a graded coassociative coproduct.
\end{lem}
\begin{pf}
Note that the above grading $deg$
is compatible with the coproduct:
$$\gamma\to\Gamma\to\gamma', \quad deg(\Gamma)=deg(\gamma)+\deg(\gamma')
=deg(\gamma\otimes\gamma').$$
For the coassociativity claim, 
see for instance \cite{CKren1}, p.12.
\end{pf}
Consider the graph homology differential (compare \cite{Kon2}, p.109, \cite{CV} p.2):
\begin{equation}\label{E:fgdiff}
d([\Gamma,l])=\sum_{e\in \Gamma^{(1)}_{int}} [\Gamma/e,l/e], 
\quad |\Gamma^{(1)}|>0.
\end{equation}
For $\Gamma\in G_{0,m}$, define $d([\Gamma,l])=0$.

To check our sign convention for the quotient orientation,
we include the following lemma.
\begin{lem}
$(H,\cdot,d)$ is dg-algebra.
\end{lem}
\begin{pf}
To prove that $d^2[\Gamma,l]=0$ (essentially as in \cite{Pen1}, p.9), 
note that in the corresponding sum:
$$\sum_{e\in\Gamma}\sum_{e'\in\Gamma,\ e'\ne e}
[(\Gamma/e)/e' , (l/e)/e'],
$$
the terms for which $e$ and $e'$ are disjoint
(the two collapsed edges do not have common vertices)
cancel in pairs, since:
$$[\Gamma/e/e', l/e/e']=-[\Gamma/e'/e, l/e'/e].$$
The other terms have an underlying graph obtained by collapsing a
two-edge subgraphs $\gamma$.
For such a fixed subgraph $\gamma$,
the coefficient $c$ of the sum of orientation classes with underlying graph 
$\Gamma/\gamma$ is the coefficient of the one-vertex graph in $G_{1,0}$: 
$d\gamma=c [\bullet,1]$.

It can be checked that it is zero for all possible 
adjacency matrices of $\gamma$.

All that is left to prove is that $d$ is a graded derivation.
It is clear that:
\begin{alignat}{1}
d([\Gamma_1,l_1][\Gamma_2,l_2])&=
\sum_{e\in\Gamma_1}[(\Gamma_1/e)\cup\Gamma_2,(l_1\vee l_2)/e]+
\sum_{e\in\Gamma_2}[\Gamma_1\cup (\Gamma_2/e), (l_1\vee l_2)/e]\notag\\
&=d([\Gamma_1,l_1])[\Gamma_2,l_2]+(-1)^{deg(\Gamma_1)}
[\Gamma_1,l_1] d([\Gamma_2,l_2]),\notag
\end{alignat}
since the collapsed edge will belong to one graph or the other,
while if $e\in\Gamma_2$:
\begin{alignat}{1}
[(\Gamma_1\cup\Gamma_2)/e,(l_1\vee l_2)/e]
&=(-1)^{deg(\Gamma_1)deg(\Gamma_2)}[(\Gamma_1\cup\Gamma_2)/e,(l_2\vee l_1)/e]\notag\\
&=(-1)^{deg(\Gamma_1)deg(\Gamma_2)}[\Gamma_1\cup(\Gamma_2/e), (l_2/e)\vee l_1]\notag\\
&=(-1)^{deg(\Gamma_1)}[\Gamma_1\cup(\Gamma_2/e),l_1\vee (l_2/e)]\notag\\
&=[\Gamma_1,l_1][\Gamma_2/e,l_2/e].\notag
\end{alignat}
\end{pf}
\begin{notation}
To simplify notation,
$[\Gamma,l]$ will be abbreviated as $\Gamma$,
the orientation class and the labeling being implicitly understood.
\end{notation}
\begin{th}\label{T:Fdgca}
$(H,d,\Delta)$ is a differential graded coalgebra.
\end{th}
\begin{pf}
All we need to prove is that $d$ is a coderivation.
Equivalently it is enough to prove that $d$ is a coderivation relative to
the reduced coproduct $\Delta_b$:
$$\Delta_b d=(d\otimes id+id\otimes d)\Delta_b.$$
In order to compare the two sides,
introduce the following correspondence.
If $\gamma\subset\Gamma$,
then denote by $\bar{\gamma}$ the corresponding
subgraph in $\Gamma/e$ ($\bar{\gamma}=\gamma/e$, if $e\in\gamma$).
Evaluating the right hand side on $\Gamma$, 
and rearranging the sum yields:
\begin{align}
RHS&=\sum_{\gamma\subset\Gamma}
d\gamma\otimes (\Gamma/\gamma)+
\sum_{\gamma\subset\Gamma}
(-1)^{deg(\gamma)}\gamma\otimes d(\Gamma/\gamma)\\
&=\sum_{\gamma\subset\Gamma}
(\sum_{e\in\gamma}\gamma/e\otimes \gamma'+
\sum_{e\in(\Gamma/\gamma}(-1)^{deg(\gamma)}\gamma\otimes (\Gamma/\gamma)/e)\notag\\
&=\sum_{e\in\Gamma}
(\sum_{e\in\gamma\subset\Gamma\to\gamma'}\gamma/e\otimes(\Gamma/\gamma)+
\sum_{e\notin\gamma\subset\Gamma\to\gamma'}\gamma\otimes((\Gamma/e)/\bar{\gamma}),
\end{align}
where a compensating sign appears in the sum for which $e\notin\gamma$,
due to the reversal of the order of taking quotients:
$\Gamma/\gamma/e$ versus $\Gamma/e/\gamma$.

Similarly:
\begin{align}
LHS&=\sum_{e\in\Gamma} \sum_{\bar{\gamma}\subset \Gamma/e}
\bar{\gamma}\otimes (\Gamma/e/\bar{\gamma})\notag\\
&=\sum_{e\in\Gamma}(
\sum_{(e/e)\in\bar{\gamma}\subset\Gamma/e}\bar{\gamma}\otimes (\Gamma/e/\bar{\gamma})+
\sum_{(e/e)\notin\bar{\gamma}\subset\Gamma/e}\bar{\gamma}\otimes ((\Gamma/e)/\bar{\gamma}),
\end{align}
where the sum was split according to whether $\bar{\gamma}$ 
contains the vertex $e/e$ (collapsed edge), or not.
With the above correspondence $\gamma\mapsto \bar{\gamma}$,
the two sides coincide,
concluding the proof.
\end{pf}
With the natural extension of the product of $H$ (abbreviated as concatenation),
to the tensor algebra $T^a(H)$,
the coproduct is an algebra morphism.
\begin{prop}
$(H,\cdot,\Delta)$ is a bialgebra.
\end{prop}
\begin{pf}
The coproduct is an algebra morphism:
$$\Delta(\Gamma_1 \Gamma_2)
=\sum_{\gamma\subset \Gamma_1\Gamma_2}\gamma\otimes (\Gamma_1 \Gamma_2)/\gamma
=\sum_{\gamma_1,\gamma_2} (\gamma_1\gamma_2)\otimes (\gamma'_1\gamma'_2)
=\Delta(\Gamma_1)\Delta(\Gamma_2).$$
where $\gamma_i=\gamma\cap\Gamma_i$ are the ``traces'' of $\gamma$ on $\Gamma_i$,
and $\gamma'_i=\Gamma_i/\gamma_i$, with $i=1,2$.
\end{pf}
Since a recursively defined antipode in a graded bialgebra comes for free:
$$S(\Gamma)=-\Gamma-\sum\limits_{\Delta_b\Gamma}S(\gamma)\gamma',$$
where the sum corresponds to the reduced coproduct $\Delta_b$,
we obtain the following result.
\begin{th}
$(H,\cdot,\Delta,S,d)$ is a dg-Hopf algebra.
\end{th}
We will be exploiting the dg-coalgebra structure on $H$,
leading to our next concern.

\section{Graph cohomology}\label{S:gc}
Returning to the main issue in this article, 
the cohomology of graphs will be derived from the dg-coalgebra structure, 
via the cobar construction 
(\cite{GLS}, p.366, \cite{M}, p.171).
%
For the reader's convenience, we will review the construction,
following \cite{T}, p.21 (see also \cite{M}, p.171).

Let $C=(H,d,\Delta,\eta,\epsilon)$ denote the coaugmented counital
dg-coalgebra of graphs from the section \ref{S:hag}.
Denote by $\bar{C}=(H/\eta(k),\bar{\Delta})$ the reduced dg-coalgebra 
corresponding to the coaugmentation ideal,
with the {\em reduced coproduct} $\bar{\Delta}$.
Note that $\bar{\Delta}=\Delta_b$ (Equation \ref{E:coprod}).

Do not shift degrees,
since $H$ already sits in the right degrees:
the elements of $\eta(k)$ are of degree -1.

Form the tensor algebra $A(\bar{C})$ and equip it with
the total differential (notation from \cite{J}, 3.4.2. p.110):
\begin{equation}\label{E:cobardif}
D=D_d+D_{\bar{\Delta}},
\end{equation} 
where $D_d$ and the ``coalgebra part'' $D_{\bar{\Delta}}$ are the graded derivations
extending $d$ and $\bar{\Delta}$ to $A(\bar{C})$.
It follows that $D$ is a derivation of $A(\bar{C})$ with $D^2=0$ (\cite{J}, p.110).
\begin{defin}\label{D:fcobar}
$F(C)=(A(\bar{C}),D)$ is the {\em cobar construction}
of the coaugmented counital dg-coalgebra $(H,d,\Delta)$ of ``Feynman graphs'' $\C{G}$.
\end{defin}
%
Consider the dual $H^*$ of $C$.
Since $C$ is of finite type,
$H^*$ is an augmented dg-algebra and 
the dual of the cobar construction $(F(C)^*,\delta)$
with dual differential $\delta=D^*$,
is isomorphic to the bar construction of $H^*$.

Taking the homology of $(F(C)^*,\delta)$ yields 
the cohomology of $H$ with coefficients in $k$.
\begin{defin}
The {\em cohomology of Feynman graphs} $\C{G}$ is:
$$H^\bullet(\C{G};k)=H_\bullet(F(C)^*,\delta),$$
where $C=(H,d,\Delta)$ is the dg-coalgebra of Feynman graphs,
and $F(C)^*$ is the dual of the cobar construction.
\end{defin}
%
We will see explicitly in the next section 
that cocycles in $F(C)^*$ determine
$L_\infty$-morphisms represented as Feynman expansions.

\section{L-infinity morphism as series over graphs}\label{S:Kfm}
L-infinity morphisms are formal diffeomorphisms between formal manifolds.
A special case of a graph expansion of an L-infinity morphism
is represented by the Kontsevich formality morphism.

The obstruction for a graded morphism between 
two DGLAs to be an L-infinity morphism will be computed,
while pointing to the physical interpretation.
As a consequence,
the coefficients of the expansion form a cocycle,
and the L-infinity morphism represented as a series over graphs
determines a cohomology class of the dg-coalgebra of graphs.

%
%
\subsection{Two DGLAs}
The context is that of \cite{Kon1}
(see also \cite{Mo}, p.79 for a slightly more general context).

Let $(A,\cdot)$ be the commutative algebra $C^\infty(X)$,
with $X=R^{d}$ for simplicity,
and $(\C{X},[,])$ a Lie subalgebra of ``vector fields'' 
of the Lie algebra of its derivations $Der(A)$.

Consider first $g_1=(T_{poly},[\ ,\ ],d=0)$,
the exterior algebra: $T_{poly}=\bigwedge^\bullet\C{X}$.
It is a DGLA with Schouten-Nijenhuis bracket and trivial differential.

The second DGLA, $g_2=(D_{poly},[,], d)$,
is the subalgebra of reduced (local) Hochschild cochains
of the Hochschild DGLA of $A$.
The cochains will be still called polydifferential operators on $A$,
since they vanish on constants.

It is assumed that $\C{X}$ has a basis $\{\partial_i\}_{i=1..d}$ over $A$,
consisting of commuting vector fields.

%
%
\subsubsection{The pre-Lie operations}
Consider coordinates $(x_1,...,x_d;y_1,...,y_d)$ in the tangent bundle $TR^d$.
In a ``supermanifold'' vain (following \cite{Kon1}, p.15),
shift the degree by one, introducing odd variables $\psi_i=y_i[1]$.
\begin{lem}
The following pre-Lie operation on vector fields induces
the Schouten-Nijenhuis Lie bracket:
$$\gamma_1\bullet \gamma_2=
\sum_{i=1}^d \frac{\partial\gamma_1}{\partial\psi_i}
\frac{\partial\gamma_2}{\partial x_i}.$$
\end{lem}
\begin{pf}
Since $[\ ,\ ]_{SN}$ is induced on the exterior algebra 
by the commutator bracket on vector fields (derivation in each variable),
while $\frac\partial{\partial\psi_i}$ and $\frac\partial{\partial x_i}$
extend as derivations on $\bigwedge^\bullet(T_{poly})$ ,
it is enough to check the claimed relation on vector fields:
$$\gamma_1=\sum_{j=1}^d \xi^j\psi_j,\quad \gamma_2=\sum_{j=1}^d\eta^j\psi_j.$$
Obviously:
$$\sum_{i=1}^d \frac{\partial\gamma_1}{\partial\psi_i}
\frac{\partial\gamma_2}{\partial x_i}=
\sum_{i=1}^d \xi^i\sum_{j=1}^d \frac{\partial\eta^j}{\partial x_i}\frac\partial{\partial x_i},$$
and then:
$$\xi\bullet\eta-\eta\bullet\xi=[\xi,\eta].$$
\end{pf}
Explicitly,
if $\gamma_1=\xi_0\wedge...\wedge\xi_k$ and $\gamma_2=\eta_0\wedge...\wedge\eta_l$,
then \cite{Mo}, p.81:
$$\gamma_1\bullet\gamma_2=\sum_{i,j=0}\epsilon(i,j) \xi_i(\eta_j)\wedge
\bar{\xi}_i\wedge\bar{\eta}_j,$$
where $\bar{\xi}_i=...\xi_{i-1}\wedge\xi_{i+1...}$ (similarly for $\bar{\eta}_j$),
$\xi_i(\eta_j)=\sum_k \xi^k (\partial_k\eta_j^l)\partial_l$,
and $\epsilon(i,j)$ is the corresponding $\pm$ sign.

Recall that the Gerstenhaber bracket corresponds to the commutator of coderivations,
under the correspondence:
$$Hom_{gr-Vect}(C,A)\cong Hom_{coalg}(C,C),$$
where $C=\oplus_{n\ge 1}A[1]^n$ is the free coassociative graded coalgebra with counit
cogenerated by the graded vector space $A[1]$ (\cite{Kon1}, p.9).
More precisely,
the graded vector space $g$ underlying the Hochschild complex:
$$C_\bullet(A;A)=Hom_{gr-Vect}(C,A),$$ 
under the above correspondence,
is the Lie algebra of coderivations of $C$.
The differential is $d=ad_{m_A}$ where $m_A$ is the associative
(and commutative) multiplication on $A$.

It is well known that the Gerstenhaber composition is a pre-Lie operation for the
Gerstenhaber bracket:
$$[\Phi_1,\Phi_2]=\Phi_1\circ\Phi_2-(-1)^{k_1k_2}\Phi_2\circ\Phi_1.$$ 

\subsection{Graph expansions and Feynman rules}\label{S:gexp}
Consider a degree -1 coderivation $U:T(g_1)\to T(g_2)$ (``pre-L-infinity morphism'')
between the two DGLAs from Section \ref{S:Kfm}.
%
It is determined by the family of skew-symmetric linear maps of degree -1
with homogeneous components (\cite{Kon1}, p.11):
$$U_n^{(k_1,...,k_n)}:g_1^{k_1}\wedge ...\wedge g_1^{k_n} \to 
g_2^{\sum k_i+1-n},
\quad n\ge 1, k_i\ge 1,$$
where $g_1^{k_i}=T^{k_i+1}(X)$ and $g_2^{m-1}=Hom(A^m,A)$.
%
\begin{rem}
Note that the degree constraint $\sum_i k_i+1-n=m-1$ requires
the number of edges of $\Gamma$, i.e. $\sum (k_i+1)$, 
to be $2n+m-2$.
\end{rem}
To expand the Taylor coefficients of the above map as a sum over graphs,
it is necessary to introduce ``states on graphs'',
beyond the labeling of vertices.
\begin{defin}\label{D:states}
Denote by $G_{n,m}^l$ the class of graphs in $\C{G}$ with
$n$ internal vertices, $m$ boundary vertices, and
$2n+m-2+l$ number of edges (\cite{Kon1}, p.23).

A {\em vertex state} on a labeled graph $(\Gamma,l)$
(not its orientation class),
is an assignment of a an element $\gamma_i\in g_1$ of degree $k_i$
for each vertex $v$ of $\Gamma$, where $i=l(v)$ and
$k_i$ is the number of outgoing edges.
The ordered sequence $(k_1,...,k_n)$ will be called the {\em signature} (type)
of the labeled graph $(\Gamma,l)$. 
\end{defin}
Fix once and for all a base $\{\partial_i\}_{i=1...d}$ of $g_1$,
of commuting vector fields.
Then any (additional) labeling of the edges defines a {\em basic state} on $\Gamma$:
$\partial_{l(e)}, e\in\Gamma^{(1)}$.
Together with a vertex state, 
a basic state defines a {\em state} on $\Gamma$, 
denoted in what follows by $\phi$.

The above pairing of a labeled graph $(\Gamma,l)$
with an element of the exterior algebra $\bigwedge(g_1)$,
factors through the projection onto its orientation class $[\Gamma,l]$
as follows.
The resulting pairing may be thought of as a 
``canonical generalized Feynman rule'' \cite{Kon1}, p.23:
\begin{equation}\label{E:U}
\C{U}:H\to Hom(\bigwedge(g_1),g_2).
\end{equation}
If $[\Gamma,l]\in G_{n,m}^l$ is a oriented graph,
and $I:E_\Gamma\to \{1,...,d\}$
is a basic state on the graph $\Gamma$,
then the only non-zero component of $\C{U}_{([\Gamma,l],I)}$ is 
\begin{equation}\label{E:gFr}
\C{U}_{([\Gamma,l],I)}(\gamma_1\wedge...\wedge\gamma_n)(f_1\otimes...\otimes f_m)
=\prod_{v\in\Gamma^{(0)}}\Phi(v),\quad \cite{Kon1}, p.23,
\end{equation}
where $\gamma_i$ are polyvector fields of degree $k_i=|Out(v_i)|$,
and $f_i\in A$.
The ``scattering amplitude'' at an interior (exterior) vertex $v$ is:
$$\Phi(v)=\Phi_{In}(v)<\gamma_{l(v)}|\Phi_{Out}(v)> \quad (\Phi(v)=\Phi_{In}(v) f_{l(v)}),$$
where the ``in'' and ``out'' ``scattering states'' are:
$$\Phi_{In}(v)=\prod_{e:\ edge\ entering\ v}\partial_{l(e)}, \quad
\Phi_{Out}(v)=\bigotimes_{e:\ edge\ out\ of\ v} dx_{l(e)}.$$
The pairing $<\ ,\ >$ between polyvector fields and forms 
extracts the corresponding coefficient
of the polyvector field attached to the vertex.

Finally define:
$$\C{U}_{[\Gamma,l]}=\sum_I \C{U}_{([\Gamma,l],I)},\quad [\Gamma,l]\in H,$$
summing over all basic states $I$ corresponding to a given vertex state.
Typical of state sums,
the result is independent of the chosen base of $g_1$.

As announced earlier,
the above ``Feynman rule'' \ref{E:gFr} is well-defined,
since the pairing of labeled graphs $(\Gamma,l)\in H_L$ with 
elements of the tensor algebra is compatible with
the action of the permutation group:
$$\C{U}':H_L\to Hom(T(g),D),\qquad 
\C{U}'_{(\Gamma,\sigma\circ l)}(\gamma)=\C{U}'_{(\Gamma,l)}(\sigma\gamma),
\quad\sigma\in\Sigma_n.$$
The pairing $\C{U}'$ will factor through the natural projections $\pi$:
$$\diagram
H_L \dto_{\pi} \rto^{\C{U}'\qquad} & Hom(T(g),D) \ar@<2pt>[d]^{\tau^*} & & T(g) \ar@<2pt>[d]^{p}\\
H \rto^{\C{U}\qquad} & Hom(\bigwedge(g),D) \ar@<2pt>[u]^{p^*} & & \bigwedge(g) \ar@<2pt>[u]^{\tau},
\enddiagram$$
since $[\Gamma,\sigma\circ l]=\epsilon(\sigma)[\Gamma,l]$,
where $\epsilon(\sigma)$ denotes the signature of the permutation $\sigma$.
Here $\tau$ denotes the usual section
$\tau(\gamma_1\wedge...\wedge\gamma_n)=\frac1{n!}\sum \epsilon(\sigma) 
\gamma_{\sigma(1)}\otimes...\otimes \gamma_{\sigma(n)}.$
\begin{rem}
In this way the above $\C{U}$ pairs labeled graphs with states, 
yielding polydifferential operators (``Feynman values'').
Note that if $\Gamma\in G_{n,m}^l$ is of type $(k_1,...,k_n)$, then
$U_\Gamma:T_{poly}^n\to D_{poly}[1-n-l]$, since:
$$U_\Gamma: g_1^{k_1}\wedge...\wedge g_n^{k_n}\to Hom(A^m,A), \quad m=2-n-l+\sum k_i.$$
If $l=0$ then $U_\Gamma$ is a pre-Lie morphism $T_{poly}[1]\to D_{poly}[1]$ (\cite{Kon1}, p.11).
\end{rem}
With the above definitions,
this canonical Feynman rule $\C{U}$ is compatible with the pre-Lie operations,
as shown next.

\begin{lem}\label{L:bullet}
If $(\Gamma',l')=(\Gamma,l)/(ij)$ is obtained by collapsing 
the {\bf interior} edge $(ij)$ of $(\Gamma,l)$, 
then:
$$\C{U}_{\Gamma}(\gamma)=\C{U}_{\Gamma'}
(\gamma_1...\hat{\gamma_i}(\gamma_i\bullet\gamma_j)...\hat{\gamma_j}...\gamma_n).$$
\end{lem}
\begin{pf}
Here the state on the collapsed interior edge (not meeting the boundary) is:
$$\phi(c)=\gamma_i\bullet\gamma_j.$$
Typical of any state-sum construction,
we only need to check the relation at the level of basic states,
after fixing a basic state outside the collapsed edge $I$, 
while summing over the basic state $k=I(ij)$ on the ``internal part of the system'',
i.e. the collapsed edge:
$$\sum_{k}\C{U}_{(\Gamma,I)}=\prod_{v\ne i,j}\Phi(v)
\sum_{k}\Phi(i)\Phi(j),\quad
\Phi(i)=\Phi_{In}(i)\gamma_i^{...k...}, \quad \Phi(j)=...\partial_k...\gamma_j^{...},$$
and:
$$\sum_k\Phi(i)\Phi(j)=\Phi(c).$$
\end{pf}
\begin{rem}
Although not entirely supported by the above pairing mechanism,
the above lemma is a generalization of Gerstenhaber $\circ_i$ composition.
If $\gamma_i\bullet\gamma_j=``\C{U}_e(\gamma_1\otimes\gamma_2)''$,
where $e\to\Gamma\to \Gamma'$,
then $\C{U}_\Gamma(\gamma_1...\gamma_n)=
\C{U}_{\Gamma'}(\gamma_1...\C{U}_e(\gamma_i\otimes\gamma_j)...\gamma_n).$
\end{rem}
The following consequence is claimed (see also \cite{Kon1}, 6.4.1.1., p.25).
\begin{cor}\label{C:bullet}
For any labeled graph $\Gamma'\in\C{G}$:
$$\sum_{e\to\Gamma\to\Gamma',\ e\in\Gamma^{(1)}_{int}}\C{U}_\Gamma(\gamma)=
\sum_{i\ne j}\C{U}_{\Gamma'}(...\wedge(\gamma_i\bullet\gamma_j)\wedge...).$$
\end{cor}
\begin{pf}
Using the previous lemma, 
the left hand side transforms as follows:
$$\sum_{e\to\Gamma\to\Gamma',\ e\in\Gamma^{(1)}_{int}}\C{U}_\Gamma(\gamma)=
\sum_{e\to\Gamma\to\Gamma', \ e\in\Gamma^{(1)}_{int}}
\C{U}_{\Gamma/(ij)}(...(\gamma_i\bullet\gamma_j)...),$$
where the non-zero terms correspond to graphs $\Gamma$ of type $(k_1,...,k_n)$,
with $k_i=deg(\gamma_i)$,
and $i,j$ are the labels of the collapsed internal edge $e$ of $\Gamma$.
Since $\Gamma/(ij)=\Gamma'$, the sum equals the right hand side.
\end{pf}
%
%
The following basic property of $\C{U}$ is expected
(``Euler-Poincare map'').
\begin{lem}\label{L:pp}
If $\Gamma_2\hookrightarrow \Gamma\twoheadrightarrow \Gamma_1$, 
where $\Gamma_1$ intersects the boundary of $\Gamma$, 
then:
$$\C{U}_\Gamma=\C{U}_{\Gamma_1}\overset{\circ}{\wedge} \C{U}_{\Gamma_2}.$$
\end{lem}
For convenience the following notation was used, with $k+l=n$:
\begin{equation}\label{E:alt}
U_k\overset{\circ}{\wedge}U_l(\gamma_1\wedge...\wedge\gamma_n)=
\frac1{(k!l!)}\sum_{\sigma\in\Sigma_n}
U_k(\gamma_{\sigma_1}\wedge...\wedge\gamma_{\sigma_k})\circ
U_l(\gamma_{\sigma_{k+1}}\wedge...\wedge\gamma_{\sigma_n}).
\end{equation}

\subsection{L-infinity morphism theorem}
The $L_\infty$-morphism obstruction is computed bellow,
for DGLAs.
\begin{th}\label{T:dgla}
Let $U:T(g_1)\to T(g_2)$ be a pre-L-infinity morphism 
represented as a Feynman expansion over the class of labeled graphs $\C{G}$:
$$U_n=\sum_{m\ge 0}\sum_{\Gamma\in G^{0}_{n,m}} W(\Gamma) \C{U}(\Gamma),$$
where $W:H\to \D{R}$ is an algebra morphism,
$\C{U}$ is defined by \ref{E:U}, 
and $G^0_{n,m}$ is the subset of graphs of $\C{G}$ with $n$ internal edges,
$m$ boundary vertices and $2n+m-2$ edges. 

(i) If $Q$ denotes the appropriate $L_\infty$-structure,
then:
\begin{equation}\label{E:Lic}
[Q,U]=(\delta W) \C{U}.
\end{equation}

(ii) $U$ is an $L_\infty$-algebra morphism iff $\delta W=0$.
\end{th}
\begin{pf}
We will prove (i), 
since (ii) becomes clear after recalling that $U$ is an $L_\infty$-morphism
iff $[Q,U]=0$ (\cite{Kel}), 
or iff condition $(F)$ holds true (\cite{Kon1}, p.24):
\begin{equation}\label{E:F}
\begin{array}{c}
\sum\limits_{i\ne j}\pm 
U_{n-1}((\gamma_i\bullet\gamma_j)\wedge\gamma_1\wedge...\wedge\gamma_n)(f_1\otimes...\otimes f_m)+\\
+\sum\limits_{k,l\ge0, k+l=n}\frac1{k!l!}\sum\limits_{\sigma\in\Sigma_n}\pm
(U_k(\gamma_{\sigma_1}\wedge...\wedge\gamma_{\sigma_k})\circ
U_l(\gamma_{\sigma_{k+1}}\wedge...\wedge\gamma_{\sigma_n}))(f_1\otimes...\otimes f_m)=0.
\end{array}
\end{equation}
Here $[Q,U]=Q_1\circ U\pm U\circ Q_2$ (see \cite{Kel}, 8,9,12).
Instead of $[Q,U]$, 
we will refer to its alternative form $(F)$.
Recall that (Equation \ref{E:cobardif}):
\begin{equation}\label{E:dW}
\delta W(\Gamma)=W(D\Gamma)=W(d\Gamma)+W(\Delta_b\Gamma),
\end{equation}
and, as stated in \cite{Kon1},
the $L_\infty$-algebra condition (F) 
corresponds to $W(d\Gamma)$ (first line - \cite{Kon1} 6.4.1.1, p.25) 
and $W(\Delta_b\Gamma)$ (second line - \cite{Kon1} 6.4.2.1., p.26).

Indeed, substitute the above Feynman expansion in the equation \ref{E:F}, 
to obtain the following component acting on $A^m$:
\begin{align}
&\sum_{\Gamma'\in G_{n-1,m}^{0}}\pm W_{\Gamma'}\sum_{i\ne j}\C{U}_{\Gamma'}
((\gamma_i\bullet\gamma_j)\wedge...) \tag{F1}\\
+&\sum\limits_{k+l=n,\ p+r=m+1}\quad
\sum_{\Gamma_1\in G_{k,p}^{0},\ \Gamma_2\in G_{l,r}^{0}}\pm W_{\Gamma_1}W_{\Gamma_2}
(\C{U}_{\Gamma_1}\overset{\circ}{\wedge}\C{U}_{\Gamma_2})(\gamma)=0,\tag{F2}
\end{align}
where $\gamma=\gamma_1\wedge...\wedge\gamma_n$, with the notation 
$\overset{\circ}{\wedge}$ from Equation \ref{E:alt}. 

Collecting the corresponding coefficients, 
the above linear combination can be written as
$\sum_{\Gamma}c_\Gamma U_\Gamma$  (\cite{Kon1}, p.25).
We claim that, for fixed $n$ and $m$, 
the only non-zero coefficients correspond to $\Gamma\in G_{n,m}^{-1}$,
and then:
\begin{equation}\label{E:c}
c_\Gamma=\delta W(\Gamma).
\end{equation}
In order to compare it with the above claim (see also Equation \ref{E:dW}),
use the above lemmas and rearrange the sums.
The first line (F1) transforms as follows:
\begin{alignat}{2}
(F1)
&=\sum_{\Gamma'\in G_{n-1,m}^0}W_{\Gamma'}
  \sum_{e\hookrightarrow\Gamma\twoheadrightarrow\Gamma',\ e\in\Gamma^{(1)}_{int}}
  \pm \C{U}_\Gamma(\gamma)\quad 
  &Corollary\ \ref{C:bullet}\quad\\
&=\sum_{\Gamma\in G_{n,m}^{-1}}
  (\sum_{e\hookrightarrow\Gamma\twoheadrightarrow\Gamma',\ e\in\Gamma^{(1)}_{int}}
  \pm W_{\Gamma/e}) \ \C{U}_\Gamma(\gamma) \quad & 
  Definition \ \ref{D:states}\\
&=\sum_{\Gamma\in G_{n,m}^{-1}}W_{d\Gamma}\ \C{U}_\Gamma(\gamma)\quad 
  &Equation\ \ref{E:fgdiff}.
\end{alignat}
The second line (F2) transforms as follows:
\begin{alignat}{2}
(F2) 
& = \sum_{k+l=n, p+r=m+1}\ \sum_{\Gamma_1\in G_{k,p}^0,\ \Gamma_2\in G_{l,r}^0}
  \pm W_{\Gamma_1}W_{\Gamma_2}
  \C{U}_{\Gamma_1}\overset{\circ}{\wedge}\C{U}_{\Gamma_2}
\quad & \\
&= \sum_{\Gamma\in G_{n,m}^{-1}}
   \sum_{\Gamma_1\hookrightarrow\Gamma\twoheadrightarrow\Gamma_2,
   \ \Gamma_1\cap\partial\Gamma\ne\emptyset}
   \pm W_{\Gamma_1}W_{\Gamma_2}\
   \C{U}_{\Gamma_1}\overset{\circ}{\wedge}\C{U}_{\Gamma_2} &\\
&= \sum_{\Gamma\in G_{n,m}^{-1}}\ 
   [\sum_{\Gamma_1\hookrightarrow\Gamma\twoheadrightarrow\Gamma_2, 
   \ \Gamma_1\cap\partial\Gamma\ne\emptyset}
  \pm W_{\Gamma_1}W_{\Gamma_2}]\quad \C{U}_{\Gamma} & Lemma\ \ref{L:pp}\\
&=\sum_{\Gamma\in G_{n,m}^{-1}} W(\Delta_b\Gamma) \ \C{U}_\Gamma. & 
\end{alignat}
The second equality is due to a resummation in the $(\Gamma_1,\Gamma_2)$-plane
``diagonally''.
Also note that,
in order to maintain the $k+l=n$ constraint,
$\Gamma_1$ must collapse to a point on the boundary,
therefore justifying the additional constraint
$\Gamma_1\cap\partial\Gamma\ne\emptyset$.

Adding the two expressions yields the terms from 
the right hand side of Equation \ref{E:dW}.
\end{pf}
\subsection{An example} 
Such a cocycle is given by Kontsevich weight function $W$,
defined in terms of the angle function $\phi$ in hyperbolic space $\C{H}$
(\cite{Kon1}, p.23):
$$\delta W(\Gamma)=0, \quad \Gamma\in G_{n,m}^{-1}.$$
In order to justify the above statement,
note first that due to dimensional reasons, $W(\Gamma)$ 
is non-zero only for $\Gamma\in G_{n,m}^0$.
Secondly, 
the codimension one boundary strata of the compactification of the configuration spaces
correspond to the coproduct (\cite{Kon1}, 5.2.1., p.22):
$$\Delta=\sum_{\gamma\to\Gamma\to \gamma'}\gamma\otimes\gamma',$$
where the sum is over all proper extensions.

Then:
$$(\delta W)(\Gamma)=W(d\Gamma+\Delta_b\Gamma)
\overset{K,I}{=}W(\Delta\Gamma)
\overset{I}{=}\int_{\partial\bar{C}_\Gamma}\omega(\Gamma)
\overset{Stokes}{=}\int_{\bar{C}_\Gamma}d\omega(\Gamma)=0,$$
where the second equality holds due to the fact that some integrals vanish 
(\cite{Kon1} p.26-27: ``bad-edge'' case and ``Type S1'' $n>2$;
\cite{pqft}, p.15-17).

\subsection{L-infinity algebras}
We claim that the above result holds for arbitrary $L_\infty-algebras$.
Moreover $L_\infty$-morphisms can be expanded over a suitable class of Feynman graphs,
and their moduli space corresponds to 
the cohomology group of the corresponding dg-coalgebra of Feynman graphs.
\begin{th}\label{T:ft}
(``Feynman-Taylor'')

Let $\C{G}$ be the class of Kontsevich graphs and 
$g_1, g_2$ two $L_\infty$-algebras as above.

In the homotopy category of $L_\infty$-algebras, 
$L_\infty$-morphisms correspond to the cohomology of the corresponding Feynman dg-coalgebra:
$\C{H}o(g_1,g_2)=H^\bullet(\C{G};k).$
\end{th}
A detailed account of the above claim is postponed to another article.

\section{Relations with perturbative QFT}\label{S:relpQFT}
The result of the previous section is important from the physical
point of view.
It is known that the Kontsevich cocycle $W$ representing the coefficients of the 
formality morphism is based on a non-linear sigma model on the disk \cite{Cat1},
and the (formality) L-infinity morphism is a partition function of a QFT.
The integrals over configuration spaces $W_\Gamma$,
are Feynman integrals for a specific propagator and interaction term in the Lagrangian.

From this point of view, 
Theorem \ref{T:ft} ``classifies'' the QFTs determined by the partition functions
corresponding to L-infinity morphisms.

The algebraic properties of Feynman graphs and Feynman integrals,
as emphasized in this article and in the Connes-Kreimer approach to renormalization,
establish an ``interface'' to a mathematical model for 
the Feynman path integral quantization based on homotopical algebra:
$$"\quad \int \C{D}\gamma\ e^{S[\gamma]}\quad "=>\sum_n\sum_\gamma U_n(\gamma).$$
The left hand side is a conceptual framework which need not be implemented using  
analytical tools (integrals, measures, etc.), 
but most likely with algebraic tools, e.g. state sum models yielding TQFTs etc..
The perturbative approach guided by a formal LHS,
and usually based on integration,
implements an expansion which is still formal, 
requiring renormalization.

The BV-formalism and $L_\infty$-formalism are two such approaches
to perturbative QFT \cite{AKSZ}.

The geometric approach of the BV-formalism implements the framework
of Feynman path integral quantization method, 
starting with a classical action, 
finding the quantum BV-action 
$S_{BV}=S_{free}+S_{int}$, 
fixing the gauge and applying a reduction technique (e.g. Faddeev-Popov)  \cite{Cat1,Cat2}, 
and finally expanding the action in perturbation series
labeled by Feynman diagrams (see \cite{Riv} for a concise introduction).

The $L_\infty$-algebra approach is a direct implementation of the RHS.
From this point of view,
to write a Lagrangian in the context of $L_\infty$-algebras \cite{Wit}, p.695,
could reasonably mean that a certain formal Lagrangian in the LHS 
is implemented in the $L_\infty$-algebra context of the RHS 
by providing a Feynman category corresponding to the interaction Lagrangian, 
and a Feynman rule correspond to the free field theory.
The Green functions would be obtained as Feynman coefficients of 
$L_\infty$-morphims expanded over a given class of Feynman graphs.
In other words, 
the partition function understood as a generating function
for the Green functions, 
may be implemented as an $L_\infty$-morphism.

At a principial level one can defend the above strategy, by claiming that the RHS is
conceptually closer to the spirit of quantum theory focusing on describing correlations.
The philosophy sketched in \cite{Irem} 
reinterprets the concept of space-time as a 
receptacle of interactions/transitions between states,
and adequately modeled by ``categories with Lagrangians'',
while the LHS comes from the traditional ``manifold approach to space-time''
forcing integrals in the sense of analysis to ``converge''.

The former philosophy can be implemented by defining a ``Feynman category''
to be essentially a ``generalized cobordism category'',
with actions as functors (see \cite{Irem}).
Cobordism categories and TQFTs, tangles, operads, PROPs, and various other graphical calculi
can be restated in terms of Generalized Cobordism Categories and their representations.



\end{document}